\documentclass[12pt]{article}
\usepackage{amsfonts}
\usepackage[dvips]{graphicx}

\begin{document}

\title{\bf Adaptive grids as parametrized scale-free networks}

\date{February 2005}

\author{\bf Gianluca Argentini \\
\normalsize gianluca.argentini@riellogroup.com \\
\textit{Advanced Computing Laboratory}\\
Information \& Communication Technology Department\\
\textit{Riello Group}, 37045 Legnago (Verona), Italy}

\maketitle

\begin{abstract}
In this paper we present a possible model of adaptive grids for numerical resolution of differential problems, using physical or geometrical properties, as viscosity or velocity gradient of a moving fluid. The relation between the values of grid step and these entities is based on the mathematical scheme offered by the model of scale-free networks, due to Bar\'abasi, so that the step can be connected to the other variables by a constitutive relation. Some examples and an application are discussed, showing that this approach can be further developed for treatment of more complex situations.
\end{abstract}

\section{A scale-free network model for adaptive grids}

We consider a differential problem of advective viscous flow in a domain $\Omega$, and let $G$ a computational grid on $\Omega$ with an initial grid step $h$. We consider a sistem of cartesian coordinates for the representation of the nodes. In 1D the grid is a linear set of {\it nodes}, in 2D is a plane surface where nodes are vertices of squares with edge length equal to $h$, in 3D is a set of adjacent cubes with squares of edge $h$ as faces (see e.g. \cite{ande}, \cite{quva}). The grid is said {\it adaptive} if, during the evolution in time, its points are clustered in regions of high flow-field gradients, so that the step $h$ can change (see e.g. \cite{ande}). If {\bf x} is the vector of the cartesian coordinates of a node and $t$ is the time variable, then for an adaptive grid we have $h=h({\bf x},t)$.\\

We consider a real positive {\it grid-field} $k$ on $G$ defined on the nodes, i.e. if $i$ is an index parameter of the $i$-th node and $M$ is the total number of nodes, we have $k=(k_i)_{1\leq i\leq M}=k({\bf x}_i,t)_{1\leq i\leq M}$. Assume that this field represents a physical or logical property associated to the flow, e.g. the number of fluid particles that interact by molecular forces with the particle on $i$-th node, or the number of possible nodes where a particle can move from time $t$ to time $t+\Delta t$ using the physical rules of the flow. As other example, in the particular matter of numerical weather prediction the topography of ground can influence heavly the computation of metheorological variables; its geometry is modelled by special grid using {\it terrain-following} coordinates (see e.g. \cite{ratto}).  Another interesting possible interpretation for $k_i$ is the number of values $({\bf u}_j)$ of the velocity field, computed at time $t$ on nodes of index $j$, that we could use in a numerical scheme for computing the value $u_i$ on $i$-th node at time $t+\Delta t$ with the desidered accuracy (see \cite{linpe}).\\
\indent From these examples we can assume that $k$ has high values in the regions of $\Omega$ where the flow is turbulent or in general not laminar; in these cases, e.g., for a good accuracy of a numerical approximation of the velocity field {\bf u} at time $t+\Delta t$, the number of values $({\bf u}_j)$, computed at time $t$, for computing the value ${\bf u}_i$ at $i$-th node can be greater then the number necessary for computing {\bf u} at a node placed, at time $t$, in a region of laminar flow (see e.g. \cite{ande}). Hence, where the grid step $h$ has small values, the grid field $k$ can have high values.\\

We apply to field $k$ the methods of a {\it scale-free network}, known in scientific literature as {\it Barab\'asi networks} (\cite{baal1}), assuming that the flow at time $t+\Delta t$ on a node of the grid depends on the flow computed at time $t$ at $i$-th node with a {\it probability function} so defined:

\begin{equation}\label{prob1}
	\Pi(k_i)= \frac{k_i}{\sum_{j=1}^{M} k_j}
\end{equation}

Using a {\it continuum approach}, Barab\'asi and Albert (\cite{baal2}) assume that $k_i$ is a continuous real variable and the rate at which it changes is proportional to $\Pi(k_i)$:

\begin{equation}\label{rate1}
	\frac{\partial k_i}{\partial t} = m\Pi(k_i)
\end{equation}

In \cite{baal2} the factor of proportionality $m$ is supposed in general constant. But if we want to describe an adaptive grid as a network changing in space and time, we suppose that the parameter $m$ is in general a function depending on the node and on time variable:

\begin{equation}\label{parameter1}
	m=(m_i)_{1\leq i\leq M}=m({\bf x}_i,t)
\end{equation}

In this case the grid becomes a scale-free network of nodes with {\it fitness}, using the terminology of Bianconi and Barab\'asi (\cite{baal2}). Also, the dependence on the grid variable {\bf x} suggests to consider the possibility of a law like (\ref{rate1}) with the gradient $\nabla k_i$ at first member.\\
The analytical dependence of $m$ from ${\bf x}_i$ and $t$ is determined by the physical and logical properties of the flow, such as viscosity, density, geometry of domain $\Omega$ and the velocity field {\bf u} itself. Hence in general we suppose that the mathematical expression of $m$ is determined by a {\it constitutive relation} such as, e.g.,

\begin{equation}\label{constitutive1}
	m=m({\bf u}, \nabla {\bf u}, \mu, t)
\end{equation}

\noindent where $\mu$ is the dynamical viscosity of the fluid.

In the next section we'll present some possible examples of choices for the analytical shape of the function $m({\bf x}_i,t)$.\\

Consider $k_i=k({\bf x}_i,t)=k_i(t)$, $m_i=m({\bf x}_i,t)=m_i(t)$. From (\ref{rate1}) and the fact that by hypothesis the field $k$ is positive, we can write

\begin{equation}
	\frac{dk_i}{k_i}=\frac{m_i}{\sum_{j=1}^{M} k_j}dt
\end{equation}

\noindent If $c=(c_i({\bf x_i}))_{1\leq i\leq M}$ is the set of initial ($t=0$) values of $k$, we have

\begin{equation}
	log(k_i)-log(c_i)=\int_{0}^{t}{\frac{m_i(s)}{\sum_{j=1}^{M} k_j(s)}ds}
\end{equation}

\noindent and therefore

\begin{equation}\label{kexp1}
	k_i = c_i \hspace{0.15cm} exp\left(\int_{0}^{t}{\frac{m_i(s)}{\sum_{j=1}^{M} k_j(s)}ds}\right)
\end{equation}

The previous equation gives the relation between the value of the grid field $k$ on the $i$-th node and the variable parameter $m$. We can write it as a field equation in this way:

\begin{equation}\label{kexp2}
	k({\bf x}_i,t) = c({\bf x}_i) \hspace{0.15cm} exp\left(\int_{0}^{t}{\frac{m({\bf x}_i,s)}{\sum_{j=1}^{M} k({\bf x}_j,s)}ds}\right)
\end{equation}

Now we state the relation between the variable grid step $h=h({\bf x},t)$ and the grid field $k$. Let ${\bf x}_i$ and ${\bf x}_j$ two {\it contiguous} nodes of $G$, e.g. $j=i+1$ in 1D, or the nodes are two vertices of the same edge of a geometrical cube in 3D. If $h_{ij}=|{\bf x}_i - {\bf x}_j|$ is the value of grid step between the nodes, according to the previous discussion we assume that $h_{ij}$ depends on the inverse of the arithmetic mean of the grid field values:

\begin{equation}\label{step1}
	h_{ij} = \frac{A}{\frac{1}{2}\left(k_i + k_j\right)}
\end{equation}

\noindent where $A$ is a real positive function. If the grid is very fine, for every ${\bf x}\in\Omega$ exists a topological neighbour $U_{{\bf x}}$ such that ${\bf x}_1$, ${\bf x}_1 \in U_{{\bf x}}$. In this case we can assume $k_i \cong k_j$ and, using (\ref{kexp2}), it is convenient to use the following definition:

\begin{equation}\label{step2}
	h({\bf x},t) = \frac{A}{k({\bf x},t)}
\end{equation}

\noindent 

\section{Examples for adaptive grids}

In all the next applications we assume $A$ constant and $\sum_{j=1}^{M} k_j = S$ with $S$ constant. We discuss some possible constitutive relation for the function $m$.\\

{\it Case 1}. Let $m_i=0$ for all nodes and for every $t$. Then, from (\ref{kexp1}), follows that $k_i=c_i \hspace{0.15cm} \forall t\geq 0$. Therefore, from (\ref{step1}), we have $h_{ij}=\frac{2A}{c_i+c_j}$, that is the grid is constant in time but it can be etherogeneous. In particular, if $c_i=c > 0 \hspace{0.15cm} \forall i$, then the grid is homogeneous. This is the case appropriate for a laminar flow, for which, in general, an adaptive grid is not strictly necessary.\\

{\it Case 2}. Let $m_i=m \neq 0$. In this case 

\begin{equation}\label{case2k}
	k_i = c_i \hspace{0.15cm} e^{\frac{mt}{S}}
\end{equation}

\noindent and

\begin{equation}\label{case2h}
	h_{ij} = \frac{2A}{c_i+c_j} \hspace{0.15cm} e^{-\frac{mt}{S}}
\end{equation}

\noindent Let $m > 0$. In this case the values of grid step at time $t + \Delta t$ are smaller than the corresponding values at time $t$, so we can compute, as previously discussed, the velocity field of a flow which has a complexity growing in time. The flow is asintotically turbulent in all points of the domain $\Omega$.\\
\noindent Let $m < 0$. In this case the values of grid step at time $t + \Delta t$ are greater than the corresponding values at time $t$. The flow is asintotically laminar in all points of the domain $\Omega$.\\
\noindent We can consider the particular case for which $c_i=c > 0 \hspace{0.15cm} \forall i$ and $m << S$. In this case 
$e^{-\frac{mt}{S}} \cong 1 - \frac{mt}{S}$, therefore the grid step $h$ is homogeneous in $\Omega$ and

\begin{equation}\label{case2hom}
	h = \frac{A}{c}\left(1 - \frac{mt}{S}\right) \hspace{1cm} \forall t\in\left[0,\frac{S}{m}\right]
\end{equation}

{\it Case 3}. According to adaptive mesh refinement methods (e.g. \cite{ranna}), grid step values can strictly depend on flow velocity gradient $\nabla {\bf u}$. Also, if the shape of domain $\Omega$ is constant in time and the viscosity $\mu$ of the fluid is a variable in space (and in time), the grid $G$ should be refined in the regions of turbulence, and in particular where $\mu$ has small values (\cite{ande}). For these reasons we can assume, as first approximation, that the parameter $m$ is constant in time and depends only on the ratio $\frac{|\nabla {\bf u}|}{\mu}$, under the hypothesis $\mu > 0$. Using a formal Taylor expansion, we write

\begin{equation}\label{mexpansion}
	m\left(\frac{|\nabla {\bf u}|}{\mu}\right) \cong m_0+m_1\frac{|\nabla {\bf u}|}{\mu}+m_2\left(\frac{|\nabla {\bf u}|}{\mu}\right)^2
\end{equation}

\noindent where $m_j$ are constants in all the grid $G$.\\
If $|\nabla{\bf u}|<<\mu$, then $m \cong m_0$ and we obtain the situations described in the previous cases. The flow becomes laminar or turbulent, depending on the sign of $m_0$.\\
In the general case, assuming $\mu=\mu({\bf x})$ constant in time, from (\ref{kexp2}) we have

\begin{equation}\label{case3exp1}
	k = c \hspace{0.15cm} exp\left(\frac{1}{\mu^2S}\int_{0}^{t}{\left[\mu^2m_0 + \mu m_1|\nabla{\bf u}| + m_2|\nabla{\bf u}|^2\right]ds}\right)
\end{equation}

\noindent Let $t > 0$. Using the definitions $\left\|\nabla{\bf u}\right\|_{L^1(0,t)}=\int_{0}^{t}{|\nabla{\bf u}(s)|ds}$, $\left\|\nabla{\bf u}\right\|_{L^2(0,t)}^2=\int_{0}^{t}{|\nabla{\bf u}(s)|^2ds}$, and remembering that the coefficients $m_j$ are constants, we can expand to first order the exponential expression and obtain

\begin{equation}\label{case3exp2}
	k = c \hspace{0.15cm}\left(1 + \frac{m_0t}{S} + \frac{m_1 \left\|\nabla{\bf u}\right\|_{L^1(0,t)}}{\mu S} + \frac{m_2 \left\|\nabla{\bf u}\right\|_{L^2(0,t)}^2}{\mu^2 S}\right)
\end{equation}

\noindent From (\ref{step2}) the grid step is

\begin{equation}\label{case3h1}
	h = \frac{\mu^2SA}{c \hspace{0.15cm}\left( \mu^2S + \mu^2 m_0t + \mu m_1 \left\|\nabla{\bf u}\right\|_{L^1(0,t)} + m_2 \left\|\nabla{\bf u}\right\|_{L^2(0,t)}^2\right)}
\end{equation}

\noindent Consider a fixed $t>0$. If $\nabla {\bf u} = {\bf 0}$, then the grid step for computing the solution at time $t+\Delta t$ becomes 

\begin{equation}\label{case3h2}
	h = \frac{SA}{c \hspace{0.15cm}\left(S + m_0t \right)}
\end{equation}

\noindent In general, if $\nabla {\bf u} \neq {\bf 0}$, for those regions in $\Omega$ where at time $t$ the viscosity $\mu$ has very small values, from (\ref{case3h1}) follows that

\begin{equation}\label{case3h3}
	\lim_{\mu \rightarrow 0} h = 0
\end{equation}

\noindent This result is coherent with usual conditions of numerical stability for differential problems where the advective terms are dominant (see e.g. \cite{quva}).

\section{Application to a diffusion-transport problem}

In this section we show how, at least in a simple differential problem where the theorical solution is known, is possible to determine the grid field $k$.\\

We consider the one-dimensional diffusion-transport problem (\cite{quva})

\begin{equation}
\left\{
\begin{array}{ll}\label{bvp1}
	$$-\mu u''(x) + b u'(x) = 0 \hspace{0.2cm} \forall x\in(0,1)$$\\
	$$u(0) = 0, \hspace{0.2cm} u(1)=1$$
\end{array}
\right.
\end{equation}

\noindent where $\mu$ and $b$ are positive constants. If $\mu << b$, case of dominant transport, the theorical solution is well approximated by the function

\begin{equation}\label{sol1}
	u(x)=e^{\left(\frac{b}{\mu}(x-1)\right)}
\end{equation}

\noindent which presents a boundary layer effect at $x=1$ (see Fig.1). 

\begin{figure}[ht]\label{expFig1}
	\begin{center}
	\includegraphics[width=7cm]{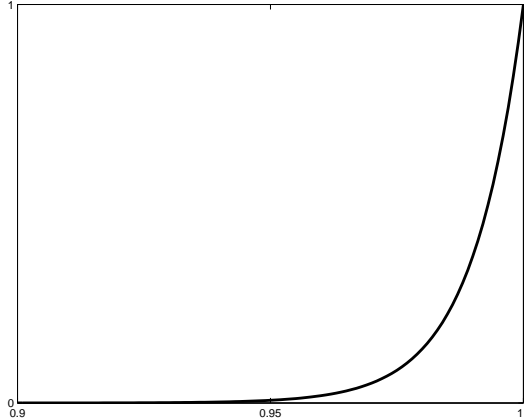}
	\caption{\small{Solution of the one-dimensional diffusion-transport problem in the case $b=100\mu$.}} 
	\end{center}
\end{figure}

\noindent For a good accuracy of a numerical resolution based on a possible numerical scheme, as Finite Differences or Finite Elements, we should impose a condition on the {\it local P\'eclet number} ${\mathbb P}{\bf e}=\frac{bh}{2\mu}$, that is ${\mathbb P}{\bf e} < 1$. For the dominant transport situation, the condition implies that a constant grid step $h$ could be very small, with large computational cost for the numerical resolution. From Fig.1 we can see that the theorical solution has significant positive values only for $x > 0.9$, therefore for the resolution of the boundary layer effect we might construct a grid with a variable step $h=h(x)$, where $h$ is decreasing on, e.g., $(0.8,1]$.\\

Because we want $h$ ever smaller when $x\rightarrow 1$, on the basis of the {\it Case 2} of the preceding section, we suppose that the entities $m$ and $\sum k_i$ are constant, $m=m_0$ and $\sum k_i = S$, but from next considerations follows that should be sufficient the ratio $\frac{m(x)}{S(x)}$ be constant. Therefore, from (\ref{kexp1}), considering the unique variable $x$ of the problem:

\begin{equation}\label{pecletk1}
 k_i=c_i \hspace{0.15cm} exp\left(\int_{0}^{x_i}\frac{m_0}{S}ds\right)=c_i \hspace{0.15cm} exp\left(\frac{m_0}{S}x_i\right)
\end{equation}

\noindent We suppose for simplicity $c_i=c \hspace{0.15cm} \forall i, 1\leq i\leq M$. Then $k_1=k(0)=c$, and from (\ref{step2}) follows that $h_1=\frac{A}{c}$, so

\begin{equation}\label{h}
	h(x)=\frac{h_1}{exp\left(\frac{m_0}{S}x\right)}
\end{equation}

\noindent Let $x_{i+1}$, $i \geq 1$, the node from which we impose the condition ${\mathbb P}{\bf e} < 1$. From definition of P\'eclet number we have 

\begin{equation}\label{d1}
	\frac{bh_1}{2\mu} < exp\left(\frac{m_0}{S}x\right) \hspace{0.5cm} \forall x\geq x_i
\end{equation}

\noindent and, with no restriction, we can impose $\frac{bh_1}{2\mu} = exp\left(\frac{m_0}{S}x_i\right)$ for the $i$-th node. In this way we obtain the required value for $\frac{m_0}{S}$:

\begin{equation}\label{m0S}
	\frac{m_0}{S}=\frac{1}{x_i}log\left(\frac{bh_1}{2\mu}\right)
\end{equation}

\noindent For our aim, the argument of the logarithmic function must to be greater then $1$, so the value of $h_1$ can be choosen so that $h_1 > \frac{2\mu}{b}$. Note that, from (\ref{bvp1}), $\frac{b}{\mu}=\frac{u''}{u'}$, so that we can consider the equation (\ref{m0S}) as a constitutive relation of kind (\ref{constitutive1}):

\begin{equation}\label{constitutive2}
	\frac{m}{S}=\frac{1}{x_i}log\left(\frac{h_1u''}{2u'}\right)
\end{equation}
	
\noindent Using (\ref{h}), the adaptive grid nodes greater than $x_i$ can be so determined (see Fig.2):

\begin{equation}\label{pecletnodes1}
	x_{j+1}=x_j + \frac{h_1}{exp\left(\frac{x_j}{x_i}log\left(\frac{bh_1}{2\mu}\right)\right)}
\end{equation}

\begin{figure}[ht]\label{expFig2}
	\begin{center}
	\includegraphics[width=6.5cm]{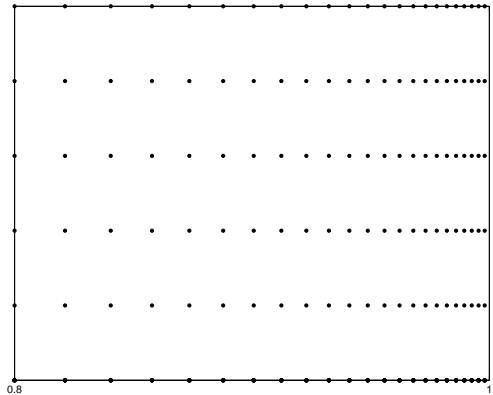}
	\caption{\small{A 2D representation of the adaptive grid for the one-dimensional diffusion-transport equation. The used values are $h_1=0.1$, $b=100\mu$, $x_i=0.8$.}}
	\end{center}
\end{figure}

\begin{figure}[ht]\label{plotFig3}
	\begin{center}
	\includegraphics[width=12cm]{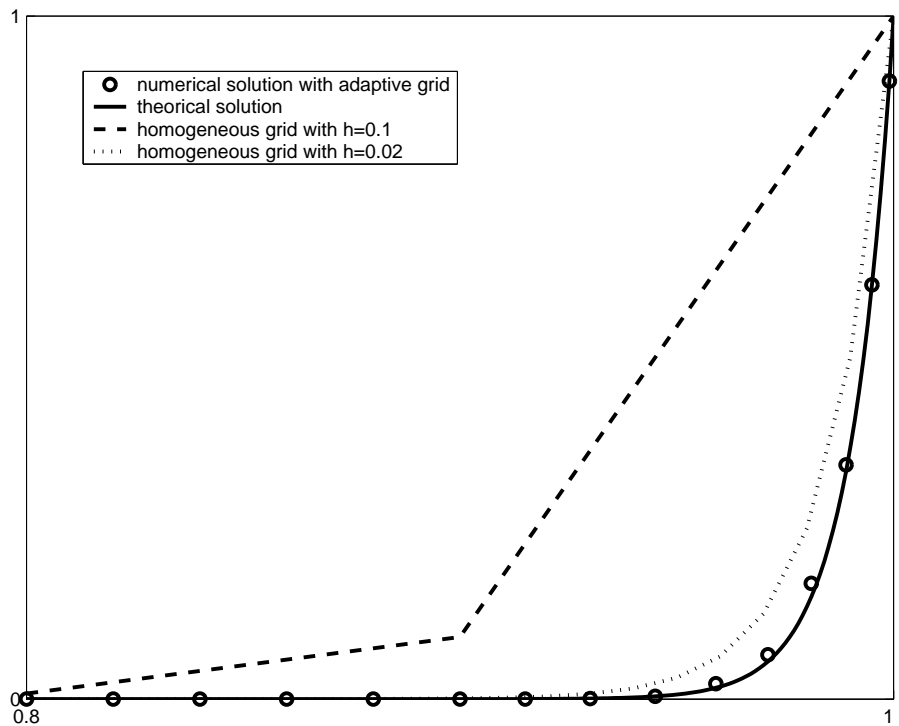}
	\caption{\small{Numerical solutions in the case $b=100\mu$, using Finite Differences method with decentrate approximation for the first derivative.}}
	\end{center}
\end{figure}

\noindent In the numerical case of Fig.\ref{plotFig3}, we have $$\lim_{x \rightarrow x_i} h(x) = 0.02, \hspace{1.5cm} \lim_{x \rightarrow 1} h = 0.0135$$ Therefore the number of nodes for this heterogeneous grid is only few greater than  the number of nodes for the homogeneous grid with $h=0.02$ and, as shown in Fig.\ref{plotFig3}, using the same numerical scheme the solution in the first case is more accurate than the solution in the latter.

\section*{Conclusion and further developments}
We have developed a model of adaptive grid based on a constituive relation between the geometrical step and other physical variables of a differential problem. The basic mathematical hypothesis is a possibile relation connecting the rate or gradient of the step and a probabilistic function which describes the physics or the geometry of the problem.\\
Further developments can regard optimization of the technical scheme and theorical treatment of the accuracy of numerical solutions obtained using this type of adaptive grids.

\end{document}